\documentclass[11pt,onecolumn,twoside]{IEEEtran} 
\usepackage{amssymb,url,amsmath,amsthm}
\usepackage{graphicx}
\usepackage{here}

\theoremstyle{definition}

\newtheorem{prop}{Proposition}

\newtheorem{theorem}{Theorem}

\def\dx{\mathrm{d}x}
\newcommand{\eqdef}{\stackrel{\mathrm{def}}{=}}

\title{Cram\'{e}r-Rao-type Bound and Stam's Inequality for Discrete Random Variables}
\author{Tomohiro Nishiyama}
\begin{document} 
\maketitle
\bibliographystyle{plain}
\begin{abstract}
The variance and the entropy power of a continuous random variable are bounded from below by the reciprocal of its Fisher information through the Cram\'{e}r-Rao bound and the Stam's inequality respectively. 
In this note, we introduce the Fisher information for discrete random variables and derive the discrete Cram\'{e}r-Rao-type bound and the discrete Stam's inequality.
\end{abstract}

\section{Introduction}
The Fisher information \cite{fisher1922mathematical} is defined for continuous random variables and plays a fundamental role in information theory and related fields.
The Fisher information for a random variable $X$ on $\mathbb{R}$ according to a probability density function $f$ is defined as
\begin{align}
\label{eq_Fisher}
I(f)\eqdef \int_\mathbb{R} \frac{f'(x)^2}{f(x)} \dx,
\end{align}
where $f'(x)$ denotes a derivative of $f(x)$ with respect to $x$.
The variance and the entropy power of a random variable $X$ are bounded from below by the reciprocal of its Fisher information through the Cram\'{e}r-Rao bound \cite{cramer2016mathematical, rao1992information} and the Stam's inequality \cite{stam1959some} respectively \cite{lutwak2005crame}.

The Cram\'{e}r-Rao bound is given by
\begin{align}
\label{CRB}
\sigma_X^2 I(f)\geq 1,
\end{align}
where $\sigma_X^2$ denotes the variance of a random variable $X$.

The Stam's inequality is given by
\begin{align}
\label{CRB} 
N(f)I(f)\geq 1,
\end{align}
where $N(f)\eqdef \frac{1}{2\pi e} \exp(2h(f))$ denotes the entropy power and $h(f)\eqdef -\int_\mathbb{R} f(x)\log f(x)\dx$ denotes the Shannon differential entropy  \cite{shannon1948mathematical}. Carlen  \cite{carlen1991superadditivity} proved the Stam's inequality is mathematically equivalent to the Gross's log-Sobolev inequality \cite{gross1975logarithmic}.  

We can also write the Fisher information as
\begin{align}
\label{eq_fisher_sqrt}
I(f)=4\int_\mathbb{R} {\biggl(\frac{d}{dx}\sqrt{f(x)}\biggr)}^2\dx \\
\label{eq_fisher_log}
I(f)=\int_\mathbb{R} f(x){\biggl(\frac{d}{dx}\log f(x)\biggr)}^2\dx.
\end{align}
Moreno, Y\'{a}nez and Dehesa introduced different discrete forms of the Fisher information  based on (\ref{eq_Fisher}), (\ref{eq_fisher_sqrt}) and (\ref{eq_fisher_log}) and they mentioned the discretization based on  (\ref{eq_fisher_sqrt}) as the most appropriate definition \cite{sanchez2009discrete}.

In this note, we introduce the diescrete Fisher information based on (\ref{eq_fisher_sqrt}) in the same way and we derive the discrete Cram\'{e}r-Rao-type bound and the discrete Stam's inequality. 

\section{Definition}
Let $\mathbb{N}_0=\mathbb{N}\cup \{0\}$.

Let $Z$ denote a discrete random variable taking values on $\mathbb{N}_0$.\\
\subsection{Set of probability mass function}
Let $p$ denote a probability mass function (pmf) of $Z$.

Let $\Omega$ be a set of probability mass functions $p$ which satisfy $\lim_{i \to \infty}p(i)=0$.
For a pmf with a finite support $[0,N-1]$, we define $p(i)=0$ for all $i\geq N$.
$\Omega$ includes many well-known distributions such as the uniform, the geometric, the Poisson, the  Bernoulli, the binomial distributions and so on.
\subsection{Discrete Fisher information (DFI)}
We introduce the discrete Fisher information (DFI) based on (\ref{eq_fisher_sqrt}). 
The DFI for $p\in\Omega$ is defined as
\begin{align}
I_d(p)\eqdef 4\sum_{i=0}^\infty \bigl(\sqrt{p(i+1)}-\sqrt{p(i)}\bigr)^2=4\sum_{i=0}^\infty D\phi(i), \nonumber
\end{align}
where $\phi(i)\eqdef \sqrt{p(i)}$ and $D$ denotes a difference operator defined as $D\phi(i)\eqdef \phi(i+1)-\phi(i)$. The advantage of the discretization based on (\ref{eq_fisher_sqrt}) is that it can be well-defined for $p(i)=0$.\\
The DFI can be also written by using the autocorrelation.
\begin{align}
I_d(p)=4(2-p(0)-2R_{\phi\phi}(1)),  \nonumber
\end{align}
where $R_{\phi\phi}(t)\eqdef \sum_{i=0}^\infty \phi(i)\phi(i+t)$ is the autocorrelation.
We can also interpret the DFI as the squared Hellinger distance between $p(i)$ and $q(i)\eqdef p(i+1)$. The squared Hellinger distance is defined as $\mathcal{H}^2(p,q)\eqdef \frac{1}{2} \sum_{i=0}^\infty\bigl(\sqrt{p(i)}- \sqrt{q(i)}\bigr)^2 $ and $I_d(p)=8\mathcal{H}^2(p,q)$ holds.
\subsection{Expected value}
The expected value for $p\in\Omega$ is defined as
\begin{align}
E[A(Z)]\eqdef \sum_{i=0}^\infty A(i) p(i), \nonumber
\end{align}
where $A$ denotes a function of a random variable $Z$.
\subsection{Entropy power}
The Shannon entropy for $p\in\Omega$ is defined as
\begin{align}
H(p)\eqdef -\sum_{i=0}^\infty p(i) \log p(i), \nonumber
\end{align}
where we define $0\log 0=0$.
The entropy power is defined as
\begin{align}
N_d(p)\eqdef \exp(2H(p)). \nonumber
\end{align}

\section{Main Results}
\subsection{New inequalities for the DFI}
\begin{theorem} (Discrete Cram\'{e}r-Rao-type bound)
\label{th_cramer-rao}
Let $p\in\Omega$ and let $I_d(p)$ be the DFI.

Let $\sigma^2\eqdef E[Z^2]-E[Z]^2$ be the variance and $\mu\eqdef E[Z]$ be the mean of a random variable $Z$ according to $p$.

Then,
\begin{align}
\label{eq_CRO}
\biggl(\sigma^2 + \frac{1}{2} - \frac{(\mu+1)^2}{2}p(0)\biggr) I_d(p)\geq \bigl(1-(\mu+1)p(0)\bigr)^2,
\end{align}
with equality if and only if $p(i)=\delta_{i0}$. $\delta_{ij}$ denotes the Kronecker delta.
\end{theorem}

When $p(0)=0$, this inequality is simplified as 
\begin{align}
\label{eq_CRO_simple}
\biggl(\sigma^2 + \frac{1}{2}\biggr) I_d(p)\geq1.
\end{align}
\begin{theorem} (Inequality for the maximum of pmf)
\label{th_max}
Let $p\in\Omega$ and let $I_d(p)$ be the DFI.

Then,
\begin{align}
\label{ineq_max}
I_d(p) > \|p\|_\infty^2 + (\|p\|_\infty-p(0))^2,
\end{align}
where $\|p\|_\infty=\max_i p(i)$.
 
Furthermore, this inequality is ``tight'' in the sense that $\alpha=1$ is the optimal constant for an inequality $  \alpha I_d(p) >\|p\|_\infty^2 + (\|p\|_\infty-p(0))^2$ which holds for all $p\in\Omega$.
\end{theorem}

\begin{prop} (Discrete Stam's inequality)
\label{prop_stam_1}
Let $p\in\Omega$ and let $I_d(p)$ be the DFI.

Then,
\begin{align}
\label{ineq_entropy_1}
N_d(p)I_d(p) > 1.
\end{align}
If there exists the optimal constant for an inequality $\beta N_d(p)I_d(p) >1$ which holds for all $p\in\Omega$, $\beta$ must be $e^{-2}\leq \beta\leq 1$.
\end{prop}

\begin{prop} (Discrete Stam-type inequality)
\label{prop_stam_2}
Let $p\in\Omega$ and let $I_d(p)$ be the DFI.

Then,
\begin{align}
\label{ineq_entropy_2}
\frac{1}{2}N_d(p)\bigl(I_d(p)+2p(0)-p(0)^2\bigr) > 1.
\end{align}
\end{prop}
When $p(0)=0$, this inequality is tighter than Proposition \ref{prop_stam_1}.

\subsection{Proofs of main results}
We show proofs of the main results.\\
\noindent\textbf{Proof of Theorem \ref{th_cramer-rao}}\\
We consider a quantity as follows.
\begin{align}
V=-\sum_{i=0}^\infty (i-\mu) Dp(i)=-\sum_{i=0}^\infty (i-\mu)(p(i+1)-p(i)) 
\end{align}
From  $\lim_{i \to \infty}p(i)=0$ and $\sum_{i=0}^\infty i p(i+1)=\sum_{i=1}^\infty (i-1) p(i)$, we have
\begin{align}
\label{eq_th1_1}
V=-\sum_{i=0}^\infty i(p(i+1)-p(i)) -\mu p(0)=\sum_{i=1}^\infty p(i) -\mu p(0)=1-(\mu+1)p(0).
\end{align}
On the other hand, we have
\begin{align}
V=-\sum_{i=0}^\infty  (i-\mu)  D[\phi(i)^2]=\sum_{i=0}^\infty  (i-\mu)  (\phi(i+1)+\phi(i))D\phi(i),
\end{align}
where we put $\phi(i)=\sqrt{p(i)}$.
Applying the Cauchy-Schwarz inequality to this equality and using the definition of the DFI  yield
\begin{align}
\label{eq_th1_2}
V^2\leq \frac{1}{4}I_d(p)\sum_{i=0}^\infty (i-\mu) ^2(\phi(i+1)+\phi(i))^2. 
\end{align}
By using $(x+y)^2\leq 2(x^2+y^2)$ and $\phi(i)=\sqrt{p(i)}$, we have
\begin{align}
\label{eq_th1_3}
\sum_{i=0}^\infty (i-\mu) ^2(\phi(i+1)+\phi(i))^2 \leq 2\sum_{i=0}^\infty (i-\mu) ^2(\phi(i+1)^2+\phi(i)^2) \\ \nonumber
=2\sum_{i=0}^\infty (i-\mu) ^2(p(i+1)+p(i))=2\sigma^2+2\sum_{i=0}^\infty (i-\mu) ^2p(i+1). 
\end{align}
Using $(i-\mu)^2=(i+1-\mu-1)^2=(i+1)^2-2(i+1)(\mu+1)+(\mu+1)^2$ yields
\begin{align}
\sum_{i=0}^\infty (i-\mu) ^2p(i+1)&=\sum_{i=0}^\infty (i+1)^2p(i+1)-2(\mu+1)\sum_{i=0}^\infty (i+1)p(i+1)+(\mu+1)^2\sum_{i=0}^\infty p(i+1) \\ \nonumber
&=E[Z^2]-2(\mu+1)\mu+(\mu+1)^2(1-p(0))=\sigma^2+1-(\mu+1)^2p(0).
\end{align}
Substituting this equality into (\ref{eq_th1_3}) and combining with  (\ref{eq_th1_2}) yields
\begin{align}
V^2\leq \biggl(\sigma^2 + \frac{1}{2} - \frac{(\mu+1)^2}{2}p(0)\biggr)I_d(p).
\end{align}
By combining this inequality with (\ref{eq_th1_1}), we obtain (\ref{eq_CRO}).

Next, we show the equality condition.
Since $(x+y)^2= 2(x^2+y^2)$ holds if and only if $x=y$, if equality holds in (\ref{eq_th1_3}), $\phi(i)=\phi(i+1)=c$ must hold for all $i$ except for $i=\mu$.
However, since $\phi(i)^2=p(i)$ satisfies $\lim_{i \to \infty}p(i)=0$, $c$ must be $0$.
Hence, if equality holds, $p(i)$ must be $\delta_{i0}$ and $\mu=0$.
By confirming the equality holds for $p(i)=\delta_{i0}$, the result follows.\\
\noindent\textbf{Proof of Theorem \ref{th_max}}\\
First, we prove the first half of the theorem.
Let $m$ be an index which satisfies $p(m)=\|p\|_\infty$.

We consider a quantity as follows.
\begin{align}
V_1=-\sum_{i=m}^\infty Dp(i)
\end{align}
From $\lim_{i \to \infty}p(i)=0$, we have
\begin{align}
\label{eq_th2_1}
V_1=p(m)=\|p\|_\infty.
\end{align}
On the other hand, we have
\begin{align}
V_1=-\sum_{i=m}^\infty(\phi(i+1)+\phi(i))D\phi(i),
\end{align}
where we put $\phi(i)=\sqrt{p(i)}$.
Applying the Cauchy-Schwarz inequality to this equality yields
\begin{align}
\label{eq_th2_2}
V_1^2\leq \sum_{i=m}^\infty |D\phi(i)|^2 \sum_{i=m}^\infty (\phi(i+1)+\phi(i))^2.
\end{align}
By using $(x+y)^2\leq 2(x^2+y^2)$, we have
\begin{align}
\label{eq_th2_3}
\sum_{i=m}^\infty (\phi(i+1)+\phi(i))^2 < 2\sum_{i=m}^\infty (\phi(i+1)^2+\phi(i)^2). 
\end{align}

Since $(x+y)^2= 2(x^2+y^2)$ holds if and only if $x=y$, if $\sum_{i=m}^\infty (\phi(i+1)+\phi(i))^2 = 2\sum_{i=m}^\infty (\phi(i+1)^2+\phi(i)^2)$ holds, $\phi(i)$ must be a constant for $i\geq m$.
Since $\phi(i)^2=p(i)$ satisfies $\lim_{i \to \infty}p(i)=0$, the constant must be $0$ and $\max_i \phi(i)=\phi(m)=0$ holds.
However,  $\max_i \phi(i)=0$ is inconsistent with $\sum_{i=0}^\infty \phi(i)^2=\sum_{i=0}^\infty p(i)=1$.
Hence, the equality doesn't hold in (\ref{eq_th2_3}). 

By using $\sum_{i=0}^\infty \phi(i)^2=1$ for (\ref{eq_th2_3}), we have
\begin{align}
\sum_{i=m}^\infty (\phi(i+1)+\phi(i))^2<4.
\end{align}
Substituting this inequality into (\ref{eq_th2_2}) and combining with (\ref{eq_th2_1}) yields
\begin{align}
\label{eq_th2_4}
4\sum_{i=m}^\infty |D\phi(i)|^2 > \|p\|_\infty^2.
\end{align}
If $m=0$, from the definition of the DFI and (\ref{eq_th2_4}), the result follows.
Then, we prove the case for $m\geq 1$ and we consider a quantity as follows.
\begin{align}
\label{eq_th2_5}
V_2=\sum_{i=0}^{m-1} Dp(i)=p(m)-p(0)=\|p\|_\infty-p(0)
\end{align}
In the same way as $V_1$, we have 
\begin{align}
\label{eq_th2_6}
V_2^2\leq \sum_{i=0}^{m-1} |D\phi(i)|^2 \sum_{i=0}^{m-1} (\phi(i+1)+\phi(i))^2 \\ \nonumber
\leq 4\sum_{i=0}^{m-1} |D\phi(i)|^2.
\end{align}
Combining  (\ref{eq_th2_5}) with (\ref{eq_th2_6}), we have
\begin{align}
\label{eq_th2_7}
4\sum_{i=0}^{m-1} |D\phi(i)|^2\geq (\|p\|_\infty-p(0))^2.
\end{align}
By taking the sum of (\ref{eq_th2_4}) and (\ref{eq_th2_7}) and using the definition of the DFI, the result follows.

Next, we prove the latter half of the theorem.
If an inequality $\alpha I_d(p) > \|p\|_\infty^2+(\|p\|_\infty-p(0))^2$ holds, $\alpha$ must satisfy $\alpha > \frac{\|p\|_\infty^2+(\|p\|_\infty-p(0))^2}{I_d(p)}$.

For the geometric distribution $p(i)=q(1-q)^i$ with $0< q \leq 1$, the DFI and the maximum of pmf are 
\begin{align}
I_d(p)=4(1-\sqrt{1-q})^2 \\ \nonumber
\|p\|_\infty=q=p(0).
\end{align}
For $q\sim 0$, from $\sqrt{1-q}= 1-\frac{q}{2}+O(q^2)$, we have
\begin{align}
\label{apr_DFI}
I_d(p)=q^2+O(q^3).
\end{align}
Hence, $\lim_{q \to +0}\frac{\|p\|_\infty^2+(\|p\|_\infty-p(0))^2}{I_d(p)}=1$ and $\alpha$ must be $\alpha \geq 1$.
Since the inequality (\ref{ineq_max}) is the case for $\alpha=1$, the result follows.\\
\noindent\textbf{Proof of Proposition \ref{prop_stam_1}}\\
First, we prove the first half of the proposition.
Since $\|p\|_\infty \geq p(i) $ and $\sum_{i=0}^\infty p(i)=1$, we have
\begin{align}
\label{eq_prop1_1}
N_d(p)= \exp\bigr(-2\sum_{i=0}^\infty p(i) \log p(i)\bigl)\geq \frac{1}{\|p\|_\infty^2}.
\end{align}
From Theorem \ref{th_max}, we have
\begin{align}
I_d(p) > \|p\|_\infty^2.
\end{align}
By combining this inequality with (\ref{eq_prop1_1}), the result follows.

Next, we prove the latter half of the theorem.
If an inequality  $\beta N_d(p)I_d(p) > 1$ holds, $\beta$ must satisfy $\beta > \frac{1}{N_d(p)I_d(p)}$.

For the geometric distribution $p(i)=q(1-q)^i$ with $0< q \leq 1$, the entropy is $H(p)=\frac{-q\log q -(1-q)\log(1-q)}{q}$.
Then, we have
\begin{align}
N_d(p)=q^{-2}(1-q)^{-\frac{2(1-q)}{q}}.
\end{align}
Combining with $\lim_{q \to +0}(1-q)^{\frac{2(1-q)}{q}}=e^{-2}$ and (\ref{apr_DFI}) yields
\begin{align}
\lim_{q \to +0}\frac{1}{N_d(p)I_d(p)}=e^{-2}.
\end{align}
Hence, $\beta$ must be $\beta \geq e^{-2}$.
Since the inequality (\ref{ineq_entropy_1}) is the case for $\beta=1$, the result follows.\\
\noindent\textbf{Proof of Proposition \ref{prop_stam_2}}\\
From $\|p\|_\infty \leq 1 $ and Theorem \ref{th_max}, we have
\begin{align}
\label{eq_prop2_1}
2\|p\|_\infty^2< I_d(p)+2\|p\|_\infty p(0)-p(0)^2\leq I_d(p)+2p(0)-p(0)^2.
\end{align}
By combining this inequality with (\ref{eq_prop1_1}), the result follows.

\section{Examples}
We show some examples of the DFI and other quantities related to the inequalities for discrete distributions.
\subsection{Discrete uniform distribution}
\begin{itemize}
\item pmf: $p(i)=\frac{1}{N}$ for $0\leq i \leq N-1$ and $p(i)=0$ for $i\geq N$. 
\item DFI: $I_d(p)=\frac{4}{N}$.
\item mean: $\mu=\frac{N-1}{2}$.
\item variance: $\sigma^2=\frac{N^2-1}{12}$.
\item maximum of pmf: $\|p\|_\infty=\frac{1}{N}$.
\item entropy power: $N_d(p)=N^2$.
\end{itemize}
\subsection{Geometric distribution}
\begin{itemize}
\item pmf: $p(i)=q(1-q)^i $ with $0< q \leq 1$.
\item DFI: $I_d(p)=4(1-\sqrt{1-q})^2$.
\item mean: $\mu=\frac{1-q}{q}$.
\item variance: $\sigma^2=\frac{1-q}{q^2}$.
\item maximum of pmf: $\|p\|_\infty=q$.
\item entropy power: $N_d(p)=q^{-2}(1-q)^{-\frac{2(1-q)}{q}}$.
\end{itemize}
\subsection{Poisson distribution}
\begin{itemize}
\item pmf: $p(i)=\frac{\lambda^i \exp(-\lambda)}{i!}$ with $\lambda > 0$.
\item DFI: $I_d(p)=4\sum_{i=0}^\infty \biggl(\sqrt{\frac{\lambda}{i+1}}-1\biggr)^2 p(i)$.
\item mean: $\mu=\lambda$.
\item variance: $\sigma^2=\lambda$.
\item maximum of pmf: $\|p\|_\infty=p(\lfloor \lambda \rfloor)$.
\item entropy power: $N_d(p)=\exp(2H(p))$ and $H(p)=\lambda(1-\log\lambda)+\exp(-\lambda)\sum_{i=0}^\infty\frac{\lambda^i \log(i!)}{i!}$.
\end{itemize}

 
\section{Conclusion}
We have introduced the discrete Fisher information (DFI) and we have shown the discrete Cram\'{e}r-Rao-type bound, the inequality for the maximum of pmf and the discrete Stam's and the Stam-type inequalities.
We have also shown the discrete Cram\'{e}r-Rao-type bound is tight and the discrete Stam's inequality is approximately tight.

It is an open question whether a tighter bound for the discrete Stam's inequality exists or not.
\bibliography{reference_Fisher}
\end{document}